\documentclass{article}
\usepackage{amsmath}
\usepackage{amsfonts}
\usepackage{amsthm}
\usepackage{amssymb}
\usepackage{hyperref}
\usepackage{mathrsfs}
\usepackage{mathtools}
\usepackage{tikz}
\usetikzlibrary{arrows.meta}
\usepackage{float}

\usepackage[numbers]{natbib}

\newcommand{\relint}{\operatorname{relint}}

\newcommand{\cubebase}{
  \tikzset{vtx/.style={circle, fill, inner sep=1.2pt}}
  \node[vtx] (e)    at ( 0,   0) {};
  \node[vtx] (v1)   at (-1.5, 1) {};
  \node[vtx] (v2)   at ( 0,   1) {};
  \node[vtx] (v3)   at ( 1.5, 1) {};
  \node[vtx] (v12)  at (-1.5, 2) {};
  \node[vtx] (v13)  at ( 0,   2) {};
  \node[vtx] (v23)  at ( 1.5, 2) {};
  \node[vtx] (v123) at ( 0,   3) {};
  \begin{scope}[gray!40, shorten >=1pt, shorten <=1pt]
    \draw (e)   -- (v1);   \draw (e)   -- (v2);   \draw (e)  -- (v3);
    \draw (v1)  -- (v12);  \draw (v1)  -- (v13);
    \draw (v2)  -- (v12);  \draw (v2)  -- (v23);
    \draw (v3)  -- (v13);  \draw (v3)  -- (v23);
    \draw (v12) -- (v123); \draw (v13) -- (v123); \draw (v23) -- (v123);
  \end{scope}}

\newcommand{\cubeflow}[3]{
  \begin{tikzpicture}[xscale=0.66667, yscale=1, baseline={(0,1.412)},
    every node/.style={font=\scriptsize}]
    \cubebase
    \foreach \a/\b/\w in {#1}
      \draw[-{Stealth[length=4pt, width=3pt]}, shorten >=1pt, shorten <=1pt]
        (\a) -- node[fill=white, inner sep=1.5pt] {\w} (\b);
    \foreach \n/\d/\p in {#2} \node[\p, inner sep=4pt] at (\n) {\d};
    \node[font=\normalsize] at (0, 3.8) {#3};
  \end{tikzpicture}}

\newcommand{\simplexf}{
  \begin{tikzpicture}[scale=0.5, baseline=(current bounding box.center),
    every node/.style={font=\scriptsize}]
    \coordinate (e1) at (0, 0);
    \coordinate (e2) at (4, 0);
    \coordinate (e3) at (2, {2*sqrt(3)});
    \coordinate (b)  at (2, {2*sqrt(3)/3});
    \draw (e1) -- (e2) -- (e3) -- cycle;
    \draw (b) -- (e1);
    \draw (b) -- (e2);
    \draw (b) -- (e3);
    \node at ({4/3}, {8*sqrt(3)/9}) {\(\mathrm{i}\)};
    \node at ({8/3}, {8*sqrt(3)/9}) {\(\mathrm{ii}\)};
    \node at (2,     {2*sqrt(3)/9}) {\(\mathrm{iii}\)};
  \end{tikzpicture}}

\newcommand{\squarebase}{
  \tikzset{vtx/.style={circle, fill, inner sep=1.2pt}}
  \node[vtx] (e)   at ( 0,  0) {};
  \node[vtx] (v1)  at (-1,  1) {};
  \node[vtx] (v2)  at ( 1,  1) {};
  \node[vtx] (v12) at ( 0,  2) {};}

\newcommand{\squarevertices}[5]{
  \begin{tikzpicture}[baseline={(0,1)}, every node/.style={font=\scriptsize}]
    \squarebase
    \node[above, inner sep=4pt] at (v12) {#1};
    \node[left, inner sep=4pt] at (v1)   {#2};
    \node[right,  inner sep=4pt] at (v2) {#3};
    \node[below, inner sep=4pt] at (e)   {#4};

    \begin{scope}[shorten >=1pt, shorten <=1pt]
      \draw (e)  -- (v1);
      \draw (e)  -- (v2);
      \draw (v1) -- (v12);
      \draw (v2) -- (v12);
    \end{scope}
    \node[font=\normalsize] at (0, 2.8) {#5};
  \end{tikzpicture}}

\newcommand{\squareedges}[5]{
  \begin{tikzpicture}[baseline={(0,1)}, every node/.style={font=\scriptsize}]
    \squarebase
    \begin{scope}[-{Stealth[length=4pt, width=3pt]}, shorten >=1pt, shorten <=1pt,
      every node/.style={font=\scriptsize, inner sep=4pt}]
      \draw (e)  -- node[below left]  {#1} (v1);
      \draw (e)  -- node[below right] {#2} (v2);
      \draw (v1) -- node[above left]  {#3} (v12);
      \draw (v2) -- node[above right] {#4} (v12);
    \end{scope}
    \node[font=\normalsize] at (0, 2.8) {#5};
  \end{tikzpicture}}

\newtheorem{conjecture}{Conjecture}[section]
\newtheorem{corollary}{Corollary}[section]
\newtheorem{lemma}{Lemma}[section]
\newtheorem{proposition}{Proposition}[section]
\newtheorem{theorem}{Theorem}[section]

\theoremstyle{definition}
\newtheorem{definition}{Definition}[section]

\theoremstyle{remark}
\newtheorem{remark}{Remark}

\title{On Kleitman's Conjecture}
\author{Jonathan Cary}
\date{\today}

\begin{document}

\maketitle
\begin{abstract}
Chv\'{a}tal conjectured that amongst the largest intersecting subfamilies of a finite subset-closed family of sets is a star.
Kleitman later strengthened Chv\'{a}tal's conjecture, defining a partial ordering on the vector space freely generated by \(2^{[n]}\) and suggesting that the vector of every maximal intersecting subfamily of \(2^{[n]}\) is bigger than a convex combination of stars.
We restate Kleitman's conjecture in terms of the cochain complex of the discrete cube, describing it as the optimization of a convex objective.
We examine the pseudoinverse of the codifferential, along with its relation to monotonicity, through which we recover the Harris-Kleitman inequality and give a generalization of a theorem of Frankl and Kupavskii on perfect matchings in superset-closed families.
We show level families to satisfy Kleitman's conjecture, providing explicit cochains.
We further show maximal intersecting families in the union of two stars and nonnegative threshold families to satisfy Kleitman's conjecture.
We close with several strengthenings of the conjecture and related open problems.
\end{abstract}

\section{Introduction}

\begin{definition}[Families of sets]
\label{families}
We say a family of sets \(\mathcal{F} \subseteq 2^{[n]}\) is intersecting if \(A \cap B \ne \emptyset\) for all \(A, B \in \mathcal{F}\), and that \(\mathcal{F}\) is a maximal intersecting family if it is intersecting and properly contained in no intersecting family.
We say \(\mathcal{F}\) is a star if \(a \in A\) for all \(A \in \mathcal{F}\) for some \(a \in [n]\).
Finally, we say \(\mathcal{F}\) is subset-closed if \(A \in \mathcal{F}\) and \(B \subseteq A\) imply \(B \in \mathcal{F}\), and superset-closed if \(A \in \mathcal{F}\) and \(A \subseteq B \subseteq [n]\) imply \(B \in \mathcal{F}\).
\end{definition}

In 1961, Erd\H{o}s, Ko, and Rado published their seminal theorem \cite{erdos1961}, that for any family \(\mathcal{F} = \binom{[n]}{k}\) such that \(2k \le n\), amongst the largest intersecting subfamilies of \(\mathcal{F}\) is a star, known as the Erd\H{o}s-Ko-Rado theorem.
In 1972, Chv\'{a}tal conjectured the following variant of the Erd\H{o}s-Ko-Rado theorem.

\begin{conjecture}
\cite{chvatal1974}
\label{chvatal}
Let \(\mathcal{F} \subseteq 2^{[n]}\) be subset-closed.
Then amongst the largest intersecting subfamilies of \(\mathcal{F}\) is a star.
\end{conjecture}

Chv\'{a}tal \cite{chvatal1974} first proved his own conjecture for all left-compressed, subset-closed families.
This result was later generalized by Wang and Wang \cite{wang1978}, then by Snevily \cite{snevily1992}, and most recently by Borg \cite{borg2011}.
Sterboul \cite{sterboul1974} proved Chv\'{a}tal's conjecture for all subset-closed families of sets of size at most 3, which was recently reproved by Czabarka, Hurlbert, and Kamat \cite{czabarka2017} via the sunflower lemma, and again by Olarte, Santos, and Spreer \cite{olarte2019} with a short combinatorial argument.
On the computational side, Eifler, Gleixner, and Pulaj \cite{eifler2022} verified the conjecture for all subset-closed families of \(2^{[7]}\) via integer programming.

Taking an alternative approach, fixing collections of intersecting families instead, Kleitman and Magnanti \cite{kleitman1974} proved the conjecture for intersecting families contained in the union of two stars.
Kleitman later proposed a conjecture in 1979, providing a general line of attack for all intersecting families, phrased in terms of the following partial ordering.

\begin{definition}[Partial ordering]
\label{kleitman-ordering}
Given \(A \subseteq [n]\), we denote by \(\tilde{A}\) the vector corresponding to \(A\) in the vector space freely generated by \(2^{[n]}\).
We define a partial ordering on this space to be that generated by the relations \(f = g + \lambda (\tilde{B} - \tilde{A})\), for \(f, g \in \mathbb{R}^{2^{[n]}}\), \(\lambda > 0\), and \(A \subset B \subseteq [n]\), and say that \(f\) is bigger than \(g\) whenever \(f\) is above \(g\) in this ordering.
\end{definition}

\begin{conjecture}
\cite{kleitman1979}
\label{kleitman}
Let \(\mathcal{F}\) be a maximal intersecting subfamily of \(2^{[n]}\).
Then \(\sum_{A \in \mathcal{F}} \tilde{A}\) is bigger than some convex combination of stars, i.e. some vector of the form \(\sum_{k \in [n]} \alpha_k \sum_{k \in A \subseteq [n]} \tilde{A}\) for some \(\alpha \in \Delta^{[n]}\).
\end{conjecture}

Relatively little progress has been made on Kleitman's conjecture, or even on Chv\'{a}tal's conjecture in the approach taken by Kleitman and Magnanti, with a recent paper by Frankl and Kupavskii \cite{frankl2023} providing another proof of Kleitman and Magnanti's result by generalizing a theorem of Berge.
Friedgut et al. \cite{friedgut2018} also recently suggested an approach for Chv\'{a}tal's conjecture in terms of correlation inequalities, with Kleitman's conjecture and several of its strengthenings being brought up in this context.

In this paper, we first develop notation for the cochain complex of the discrete cube in terms of a Fourier-Walsh basis.
The pseudoinverse of the codifferential is discussed, and we recover the Harris-Kleitman inequality and give a generalization of a theorem of Frankl and Kupavskii.
We then introduce polarity as a reformulation of Kleitman's conjecture, giving a basic analysis of the property using its convex dual, and we establish polarity for maximal intersecting families contained in the union of two stars and nonnegative threshold families.

\section{Preliminaries}

\subsection{Real-Valued Boolean Functions}

We first give preliminary definitions and propositions for real-valued Boolean functions, following those given in \cite{odonnell2014}.
Throughout, we identify subsets of \([n]\) with elements of \(\{-1, 1\}^n\) in the natural way.

\begin{definition}[Inner product and Walsh functions]
\label{walsh}
We equip the space of real-valued Boolean functions \(f \colon 2^{[n]} \to \mathbb{R}\) with the inner product
\[
    \langle f, g \rangle = \frac{1}{2^n} \sum_{A \subseteq [n]} f(A) g(A).
\]
For each subset \(A \subseteq [n]\), the Walsh function \(x^A \colon 2^{[n]} \to \mathbb{R}\) is defined by
\[
    x^A(B) = (-1)^{|A \setminus B|}, \quad \text{for all } B \subseteq [n].
\]
\end{definition}

In the case of singletons \(\{a\} \subseteq [n]\), we write \(x_a\) to denote \(x^{\{a\}}\).

The Walsh functions form an orthonormal basis with respect to the inner product above.
Consequently, any real-valued Boolean function \(f \colon 2^{[n]} \to \mathbb{R}\) admits a unique expansion
\[
    f = \sum_{A \subseteq [n]} \hat{f}(A) x^A.
\]

\begin{definition}[Fourier-Walsh transform]
\label{fourier-walsh}
The coefficients \(\hat{f}(A) = \langle f, x^A \rangle\) appearing in the expansion above are called the Fourier-Walsh coefficients of \(f\), and the map \(\hat{f} \colon 2^{[n]} \to \mathbb{R}\) is called the Fourier-Walsh transform of \(f\).
\end{definition}

\begin{definition}[Partial derivative and restriction]
\label{derivative-restriction}
The partial derivative with respect to coordinate \(k \in [n]\) is defined by
\[
    \frac{\partial}{\partial x_k} f \colon 2^{[n] \setminus \{k\}} \to \mathbb{R}, \quad
    \left( \frac{\partial}{\partial x_k} f \right)(A) = \frac{f(A \cup \{k\}) - f(A)}{2},
\]
and the restriction to \(x_k = a\) for \(k \in [n]\) by
\[
    f \big|_{x_k = a} \colon 2^{[n] \setminus \{k\}} \to \mathbb{R}, \quad
    \left( f \big|_{x_k = a} \right)(A) = \frac{(1 - a)f(A) + (1 + a)f(A \cup \{k\})}{2}.
\]
\end{definition}

\begin{definition}[Convolution]
\label{convolution}
The convolution of real-valued Boolean functions is defined by
\[
    (f * g)(A) = \frac{1}{2^n} \sum_{B \odot C = A} f(B) g(C), \text{ for all } A \subseteq [n].
\]
\end{definition}

\begin{definition}[Monotonicity, antipodality, and intersection]
\label{monotone-antipodal}
Finally, for a real-valued Boolean function \(f \colon 2^{[n]} \to \mathbb{R}\), we say \(f\) is monotone if \(A \subseteq B\) implies \(f(A) \le f(B)\), and antipodal if \(f(A) = -f([n] \setminus A)\) for all \(A \subseteq [n]\). In the case of \(f \colon 2^{[n]} \to \{-1, 1\}\), we say \(f\) is intersecting if \(f^{-1}(1)\) is intersecting.
\end{definition}

\begin{proposition}
The partial derivative satisfies \[
    \frac{\partial}{\partial x_k} f = \sum_{k \in A \subseteq [n]} \hat{f}(A) x^{A \setminus \{k\}}.
\]
\end{proposition}

\begin{proposition}
The restriction satisfies \[
    f \big|_{\{ x_k = a_k \}_{k \in B}} = \sum_{A \subseteq [n]} \hat{f}(A) \left( \prod_{k \in A \cap B} a_k \right) x^{A \setminus B}.
\]
\end{proposition}

\begin{proposition}
The convolution satisfies \[
    f * g = \sum_{A \subseteq [n]} \hat{f}(A) \hat{g}(A) x^A.
\]
\end{proposition}

\subsection{Cochains on the Discrete Cube}

We develop cochain machinery specializing the cubical cohomology framework of \cite{kaczynski2004,kaczynski2013}, adapting it to a Fourier-Walsh basis for real-valued Boolean functions.

\begin{definition}[The discrete cubical complex]
\label{cubical-complex}
We define \(3^{[n]}\) to be the discrete cubical complex of dimension \(n\), consisting of the set of all products \(\prod_{i \in [n]} I_i\), where each \(I_i\) is either the closed interval \([-1, 1]\) or one of its boundary points.
We will identify elements \(\sigma \in 3^{[n]}\) with tuples \((A, B) \in 2^{[n]}\) such that \(A \cap B = \emptyset\), where \(A\) is the set of non-degenerate axes of \(\sigma\) and \(B\) corresponds to the coordinates of the degenerate axes of \(\sigma\) in the natural way.

Let \(\binom{[n]}{k} \times 2^{[n] - k}\) for \(k \in [n]\) denote the set of all \(k\)-faces of \(3^{[n]}\), that is, all such products in which exactly \(k\) coordinates are non-degenerate intervals, and for \(A \subseteq [n]\), we write \(\{A\} \times 2^{[n] \setminus A}\) for the set of all \(\lvert A \rvert\)-faces of \(3^{[n]}\) that are non-degenerate precisely along the axes indexed by \(A\).

The face maps \(\delta_a^{\varepsilon}\) map \(k\)-cubes of \(3^{[n]}\) to their \((k - 1)\)-faces by setting the \(a\)th coordinate to \(\varepsilon \in \{-1, 1\}\), where \(a\) must be non-degenerate.

The cubical product \(\diamond\) on \(3^{[n]}\) is defined such that for any \(f, g \in 3^{[n]}\) intersecting at a single point, \(f \diamond g\) is the unique face such that each interval is the common degenerate interval if \(f\) and \(g\) agree there, and otherwise is the unique non-degenerate interval from either \(f\) or \(g\).
\end{definition}

\begin{definition}[Chain complex and boundary operator]
\label{chain-complex}
We define the chain complex \(C_{\bullet}(3^{[n]})\) to be the graded vector space freely generated by the elements of \(3^{[n]}\).
For each \(k \in [n]\), the space of \(k\)-chains, denoted \(C_k(3^{[n]})\), is the subspace generated by \(\binom{[n]}{k} \times 2^{[n] - k}\).

The boundary operator \(\partial_k \colon C_k(3^{[n]}) \to C_{k-1}(3^{[n]})\) is defined by \[
    \partial_k(\tilde{\sigma}) = \sum_{a \in A} (-1)^{\varepsilon_A(a)} \left( \widetilde{\delta_a^{1}(\sigma)} - \widetilde{\delta_a^{-1}(\sigma)} \right),
\] where \(A\) is the set of non-degenerate axes of \(\sigma\) and \(\varepsilon_A(a) = \lvert \{ b \in A \mid b < a \} \rvert\), which we extend additively to sets by \(\varepsilon_A(B) = \sum_{b \in B} \varepsilon_A(b)\).
\end{definition}

\begin{definition}[Cochain complex, pairing, and inner product]
\label{cochain-complex}
The cochain complex \(C^{\bullet}(3^{[n]})\) is defined as the graded dual of the chain complex:
\[
    C^k(3^{[n]}) = \mathrm{Hom}(C_k(3^{[n]}), \mathbb{R}),
\]
the space of real-valued linear functionals on \(k\)-chains.

The natural pairing between cochains and chains is defined by
\[
    \langle \phi, \tilde{\sigma} \rangle = \frac{1}{2^{n - k}} \phi(\tilde{\sigma})
    \quad \text{for } \phi \in C^k(3^{[n]}),\ \sigma \in \binom{[n]}{k} \times 2^{[n] - k},
\]
and we equip \(C^k(3^{[n]})\) with the inner product
\[
    \langle \phi, \psi \rangle = \frac{1}{2^{n - k}} \sum_{\sigma \in \binom{[n]}{k} \times 2^{[n] - k}} \phi(\tilde{\sigma}) \psi(\tilde{\sigma}),
\]
for all \(\phi, \psi \in C^k(3^{[n]})\).
\end{definition}

\begin{definition}[Exterior derivative and codifferential]
\label{exterior-derivative}
The exterior derivative \(d \colon C^k(3^{[n]}) \to C^{k+1}(3^{[n]})\) is defined as the adjoint of the boundary operator:
\[
    (df)(\sigma) = \frac{1}{2} f(\partial \sigma).
\]

The codifferential \(\delta \colon C^{k+1}(3^{[n]}) \to C^k(3^{[n]})\) is defined as the adjoint of \(d\) with respect to the inner product:
\[
    \langle d\phi, \psi \rangle = \langle \phi, \delta \psi \rangle \quad \text{for all } \phi \in C^k(3^{[n]}),\ \psi \in C^{k+1}(3^{[n]}).
\]
\end{definition}

\begin{definition}[Exterior product]
\label{exterior-product}
We define the exterior product on chains and on cochains by \begin{align*}
    \tilde{\sigma} \wedge \tilde{\tau} & = (-1)^{\varepsilon_B(A)} \frac{1}{2^{a+b}} \, \widetilde{\sigma \diamond \tau} \quad \text{for } \sigma \in \binom{[n]}{a} \times 2^{[n]\setminus a}, \tau \in \binom{[n]}{b} \times 2^{[n]\setminus b}, \\
    (\phi \wedge \psi)(\tilde{\sigma}) & = \frac{1}{2^{a+b}} \sum_{\alpha \diamond \beta = \sigma} (-1)^{\varepsilon_B(A)} \phi(\tilde{\alpha}) \psi(\tilde{\beta}) \quad \text{for } \phi \in C^a(3^{[n]}), \psi \in C^b(3^{[n]}),
\end{align*}
where in each case \(A\) and \(B\) denote the non-degenerate axes of the first and second factor, so that \(\varepsilon_B(A) = \lvert \{ (c, d) \in A \times B \colon d < c \} \rvert\).
\end{definition}

\begin{definition}[Basis chains and cochains]
\label{basis-cochains}
For each subset \(A \subseteq [n]\), we define the chain \(\frac{\partial}{\partial x^A} \in C_{\lvert A \rvert}(3^{[n]})\) and cochain \(dx^A \in C^{\lvert A \rvert}(3^{[n]})\) by
\begin{align*}
    \frac{\partial}{\partial x^A} & = \sum_{\sigma \in \{A\} \times 2^{[n] \setminus A}} \tilde{\sigma}, \\
    dx^A(\tilde{\sigma})          & = \begin{cases}
        1, & \text{if } \sigma \in \{A\} \times 2^{[n] \setminus A}, \\
        0, & \text{otherwise}.
    \end{cases}
\end{align*}
\end{definition}

In the case of singletons \(\{a\} \subseteq [n]\), we write \(\frac{\partial}{\partial x_a}\) and \(dx_a\) to denote \(\frac{\partial}{\partial x^{\{a\}}}\) and \(dx^{\{a\}}\) respectively.

Every cochain \(\phi \in C^{\bullet}(3^{[n]})\) then admits a multilinear representation of the form
\[
    \phi = \sum_{A \subseteq [n]} \sum_{B \subseteq [n] \setminus A} \widehat{\phi_A}(B) x^B dx^A,
\]
where \(\phi_A \colon 2^{[n] \setminus A} \to \mathbb{R}\) denotes the real-valued Boolean function corresponding to the \(dx^A\)-component of \(\phi\).
Here again, we write \(\phi_a\) in place of \(\phi_{\{a\}}\) for singletons \(\{a\} \subseteq [n]\).

\begin{definition}[Further operators]
Finally, we define the gradient by
\[
    \nabla f = \sum_{k \in [n]} \left( \frac{\partial}{\partial x_k} f \right) \frac{\partial}{\partial x_k} \quad \text{for } f \colon 2^{[n]} \to \mathbb{R},
\]
the interior product \(\iota\) by
\begin{align*}
    \iota_X \phi & = \langle \phi, X \rangle \quad \text{for } \phi \in C^1(3^{[n]}), X \in C_1(3^{[n]}), \\
    \langle \iota_X \phi, Y \rangle & = \langle \phi, X \wedge Y \rangle \quad \text{for } \phi \in C^k(3^{[n]}), X \in C_1(3^{[n]}), Y \in C_{k - 1}(3^{[n]}),
\end{align*}
the convolution of cochains by
\[
    \phi * \psi = \sum_{A \subseteq [n]} (\phi_A * \psi_A) dx^A \quad \text{for } \phi, \psi \in C^\bullet(3^{[n]}),
\]
and the Laplace-de Rham operator \(\Delta\) by
\[
    \Delta = d \delta + \delta d.
\]
\end{definition}

\begin{proposition}
The exterior product satisfies \[
    \phi \wedge \psi = \sum_{A \cap B = \emptyset} (-1)^{\varepsilon_B(A)} \left. \phi_A \psi_B \right|_{\{ x_a = 0 \}_{a \in A \cup B}} dx^{A \cup B}.
\]
\end{proposition}

\begin{proposition}
The interior product satisfies \[
    \iota_{X} \phi = \sum_{A \subseteq B \subseteq [n]} (-1)^{\varepsilon_{B \setminus A}(A)} \left. X_{A} \phi_{B} \right|_{\{ x_c = 0 \}_{c \in B \setminus A}} dx^{B \setminus A}.
\]
\end{proposition}

\begin{proposition}
The exterior derivative \(d\phi\) satisfies \begin{align*}
    d\phi & = \sum_{A \subseteq [n]} \sum_{a \in [n] \setminus A} \frac{\partial}{\partial x_a} \phi_A \, dx_a \wedge dx^A \\
          & = \sum_{A \subseteq [n]} \sum_{B \subseteq [n] \setminus A} \sum_{a \in B} (-1)^{\varepsilon_A(a)} \widehat{\phi_A}(B) \, x^{B \setminus \{a\}} dx^{A \cup \{a\}}.
\end{align*}
\end{proposition}

\begin{proposition}
The codifferential \(\delta \phi\) satisfies \begin{align*}
    \delta \phi & = \sum_{A \subseteq [n]} \sum_{a \in A} x_a \iota_{\frac{\partial}{\partial x_a}} \phi_A dx^A \\
                & = \sum_{A \subseteq [n]} \sum_{B \subseteq [n] \setminus A} \sum_{a \in A} (-1)^{\varepsilon_A(a)} \widehat{\phi_A}(B) x^{B \cup \{a\}} dx^{A \setminus \{a\}}.
\end{align*}
\end{proposition}

\begin{theorem}
\label{diagonalization}
The Laplace-de Rham operator satisfies \[
    \Delta \phi = \sum_{A \subseteq [n]} \sum_{B \subseteq [n] \setminus A} (\lvert A \rvert + \lvert B \rvert) \widehat{\phi_A}(B) x^B dx^A.
\]
\end{theorem}

\begin{proof}
By linearity it suffices to consider \(x^B dx^A\) for \(A \sqcup B \subseteq [n]\).
Noting that \begin{align*}
    d (x^B dx^A)      & = \sum_{b \in B} (-1)^{\varepsilon_A(b)} x^{B \setminus \{b\}} dx^{A \cup \{b\}} \text{ and} \\
    \delta (x^B dx^A) & = \sum_{a \in A} (-1)^{\varepsilon_A(a)} x^{B \cup \{a\}} dx^{A \setminus \{a\}},
\end{align*}
we see that \(d\) moves a coordinate from \(B\) to \(A\) and \(\delta\) from \(A\) to \(B\).
Expanding \(\Delta = d \delta + \delta d\) we thus have terms of the form \(x^B dx^A\) and \(x^{B \triangle \{a, b\}} dx^{A \triangle \{a, b\}}\) for \(a \in A\) and \(b \in B\).
We proceed by cases.
\begin{enumerate}
    \item Each \(x^B dx^A\) term arises once for each \(a \in A\) from \(d \delta\) and once for each \(b \in B\) from \(\delta d\).
    We have \(\varepsilon_{A \setminus \{a\}}(a) = \varepsilon_A(a)\) and \(\varepsilon_{A \cup \{b\}}(b) = \varepsilon_A(b)\), so the two signs in each such term agree.
    We thus have \(\lvert A \rvert + \lvert B \rvert\) terms of \(x^B dx^A\).
    \item Each \(x^{B \triangle \{a, b\}} dx^{A \triangle \{a, b\}}\) term for \(a \in A\) and \(b \in B\) arises exactly twice: once from \(d \delta\), with sign exponent \(\varepsilon_A(a) + \varepsilon_{A \setminus \{a\}}(b)\), and once from \(\delta d\), with sign exponent \(\varepsilon_A(b) + \varepsilon_{A \cup \{b\}}(a)\).
    Removing \(a\) from \(A\) lowers \(\varepsilon_A(b)\) by one exactly when \(a < b\), while adding \(b\) to \(A\) raises \(\varepsilon_A(a)\) by one exactly when \(b < a\).
    Precisely one of these occurs, so the two exponents differ in parity and the pair cancels.
\end{enumerate}
\end{proof}

\begin{corollary}
Every closed cochain is exact, and every coclosed cochain is coexact, up to a constant \(0\)-cochain.
\end{corollary}

\begin{proof}
By Theorem~\ref{diagonalization}, \(\Delta\) is invertible up to a constant \(0\)-cochain and commutes with \(d\) and \(\delta\).
Thus for closed \(\phi \in C^\bullet(3^{[n]})\) with \(\widehat{\phi_\emptyset}(\emptyset) = 0\), we have \[
    \phi = \Delta \Delta^{-1} \phi = (d \delta + \delta d) \Delta^{-1} \phi = d (\delta \Delta^{-1} \phi) + \delta \Delta^{-1} (d \phi) = d (\delta \Delta^{-1} \phi).
\]
The coclosed case follows similarly.
\end{proof}

\begin{theorem}
\label{leibniz}
For all \(f \colon 2^{[n]} \to \mathbb{R}\) and \(\phi \in C^k(3^{[n]})\), \[
    \delta (f \wedge \phi) = f \wedge \delta \phi - \iota_{\nabla f} \phi.
\]
\end{theorem}

\begin{proof}
Both sides are bilinear in \(f\) and \(\phi\), so we may take \(f = x^A\) and \(\phi = x^B dx^C\) for \(A, C \subseteq [n]\) and \(B \subseteq [n] \setminus C\).
Here \(x^A \wedge x^B dx^C = x^{A \triangle B} \big|_{\{ x_c = 0 \}_{c \in C}} dx^C\), which vanishes unless \(A \cap C = \emptyset\), while the \(k\)th term of \(\iota_{\nabla x^A} (x^B dx^C)\) survives precisely when \(A \cap C = \{k\}\).
We proceed by cases on \(\lvert A \cap C \rvert\).
\begin{enumerate}
    \item For \(\lvert A \cap C \rvert = 0\), the interior product vanishes, and we have \begin{align*}
        \delta (x^A \wedge x^B dx^C) & = \delta (x^{A \triangle B} dx^C) \\
                                     & = x^A \wedge \delta (x^B dx^C) \\
                                     & = x^A \wedge \delta (x^B dx^C) - \iota_{\nabla x^A} (x^B dx^C).
    \end{align*}
    \item For \(\lvert A \cap C \rvert = 1\), say \(A \cap C = \{k\}\), the wedge vanishes, and only the \(c = k\) term of \(\delta (x^B dx^C)\) survives wedging with \(x^A\), so \begin{align*}
        \delta (x^A \wedge x^B dx^C) & = 0 \\
                                     & = \left( (-1)^{\varepsilon_C(k)} - (-1)^{\varepsilon_C(k)} \right) x^{A \setminus \{k\}} \wedge x^B dx^{C \setminus \{k\}} \\
                                     & = x^A \wedge \delta (x^B dx^C) - \iota_{\nabla x^A} (x^B dx^C).
    \end{align*}
    \item For \(\lvert A \cap C \rvert \ge 2\), both wedges vanish, since \(A\) still meets \(C \setminus \{c\}\) for each \(c \in C\), and no \(k\) has \(A \cap C = \{k\}\), so \[
        \delta (x^A \wedge x^B dx^C) = x^A \wedge \delta (x^B dx^C) = \iota_{\nabla x^A} (x^B dx^C) = 0
    \] and thus \[
        \delta (x^A \wedge x^B dx^C) = x^A \wedge \delta (x^B dx^C) - \iota_{\nabla x^A} (x^B dx^C).
    \]
\end{enumerate}
\end{proof}

\begin{theorem}
\label{pseudoinverse}
The pseudoinverse of \(\delta \colon C^1(3^{[n]}) \to C^0(3^{[n]})\) satisfies \[
    \delta^\dagger f = \sum_{a \in [n]} \sum_{A \subseteq [n] \setminus \{a\}} \frac{1}{\lvert A \rvert + 1} \hat{f}(A \cup \{a\}) x^{A} dx_a.
\]
Equivalently,
\[
    \delta^\dagger f = (\delta^\dagger \delta_e) * df,
\]
where \(\delta_e = \sum_{A \subseteq [n]} x^A\) is the identity with respect to convolution, and
\[
    \delta^\dagger \delta_e = \sum_{a \in [n]} \sum_{A \subseteq [n] \setminus \{a\}} \frac{1}{\lvert A \rvert + 1} x^{A} dx_a = \sum_{A \subseteq [n]} \frac{1}{(n - \lvert A \rvert)!} \sum_{\phi \in \Phi(A)} \phi,
\]
where \(\Phi(A)\) denotes the set of all \(1\)-cochains of directed paths in the discrete cube from \(A\) to \([n]\).
\end{theorem}

\begin{proof}
The codifferential satisfies \(\delta(x^A dx_a) = x^{A \cup \{a\}}\).
Break \(C^1(3^{[n]})\) into the spans \(V_B\) of \(\{ x^{B \setminus \{b\}} dx_b : b \in B \}\) and \(C^0(3^{[n]})\) into the spans of \(x^B\), for each \(B \subseteq [n]\); these decompositions are orthogonal, and \(\delta\) carries each of the \(\lvert B \rvert\) basis vectors of \(V_B\) to \(x^B\).
The claimed \(\delta^\dagger f\) has all \(\lvert B \rvert\) coordinates in \(V_B\) equal to \(\hat{f}(B) / \lvert B \rvert\), so
\begin{itemize}
    \item \(\delta \delta^\dagger f = \sum_{\emptyset \subset B \subseteq [n]} \hat{f}(B) x^B = f - \hat{f}(\emptyset)\), the projection of \(f\) onto \(\operatorname{range}(\delta)\), and
    \item \(\delta^\dagger f \perp \ker(\delta)\), since \(\ker(\delta) \cap V_B\) consists of the coordinate vectors summing to zero.
\end{itemize}
These are precisely the conditions characterizing \(\delta^\dagger f\) as the least-norm minimizer of \(\lVert \delta \phi - f \rVert\).

Finally, to see the equivalence of the two expressions of \(\delta^\dagger \delta_e\), note that both are closed with equal codifferentials.
\end{proof}

The three cochains \(f\), \(df\), and \(\delta^\dagger f\) are recorded on the discrete square in Fig.~\ref{square-figure}.

\begin{figure}[!b]
\centering
\setlength{\tabcolsep}{6pt}
\begin{tabular}{ccc}
    \squarevertices{\(a\)}{\(b\)}{\(c\)}{\(d\)}{\(f\)} &
    \squareedges
        {\(\frac{b - d}{2}\)}{\(\frac{c - d}{2}\)}
        {\(\frac{a - b}{2}\)}{\(\frac{a - c}{2}\)}{\(df\)} &
    \squareedges
        {\(\frac{a + 3b - c - 3d}{8}\)}
        {\(\frac{a - b + 3c - 3d}{8}\)}
        {\(\frac{3a - 3b + c - d}{8}\)}
        {\(\frac{3a + b - 3c - d}{8}\)}{\(\delta^\dagger f\)}
\end{tabular}
\caption{\(f\), \(df\), and \(\delta^\dagger f\) on the discrete square.}
\label{square-figure}
\end{figure}

\begin{corollary}
\label{pseudoinverse-nonnegative}
Let \(f \colon 2^{[n]} \to \mathbb{R}\) be monotone.
Then \(\delta^\dagger f\) is nonnegative.
\end{corollary}

\begin{proof}
Here \(df\) is nonnegative by monotonicity, and \(\delta^\dagger \delta_e\) is nonnegative by the monotone path expression of Theorem~\ref{pseudoinverse}.
As a convolution of nonnegative functions is nonnegative, so too is \(\delta^\dagger f = (\delta^\dagger \delta_e) * df\).
\end{proof}

From Corollary~\ref{pseudoinverse-nonnegative}, the Harris-Kleitman inequality and a generalization of a theorem of Frankl and Kupavskii follow.

\begin{theorem}[Harris-Kleitman inequality {\cite{harris1960,kleitman1966}}]
\label{harris-kleitman}
Let \(\mathcal{F}, \mathcal{G} \subseteq 2^{[n]}\) be superset-closed.
Then \(\frac{\lvert \mathcal{F} \cap \mathcal{G} \rvert}{2^n} \ge \frac{\lvert \mathcal{F} \rvert}{2^n} \frac{\lvert \mathcal{G} \rvert}{2^n}\), with equality if and only if \(\mathcal{F}\) and \(\mathcal{G}\) are determined by disjoint sets of elements. 
\end{theorem}

\begin{proof}
Let \(f = \sum_{F \in \mathcal{F}} \tilde{F}\) and \(g = \sum_{G \in \mathcal{G}} \tilde{G}\).
Then \begin{align*}
    \frac{\lvert \mathcal{F} \cap \mathcal{G} \rvert}{2^n} & = \langle f, g \rangle \\
                                                           & = \langle df, \delta^\dagger g \rangle + \hat{f}(\emptyset)\hat{g}(\emptyset) \\
                                                           & \ge \hat{f}(\emptyset)\hat{g}(\emptyset) \\
                                                           & = \frac{\lvert \mathcal{F} \rvert}{2^n} \frac{\lvert \mathcal{G} \rvert}{2^n},
\end{align*}
where \(\langle df, \delta^\dagger g \rangle \ge 0\) since \(df \ge 0\) by monotonicity and \(\delta^\dagger g \ge 0\) by Corollary~\ref{pseudoinverse-nonnegative}.
Moreover, we have \(\langle df, \delta^\dagger g \rangle = 0\) if and only if for each \(k \in [n]\), at least one of the partials, \(\frac{\partial}{\partial x_k} f\) and \(\frac{\partial}{\partial x_k} g\), is \(0\).
\end{proof}

\begin{remark}
In fact, any nonnegative \(\phi \in C^1(3^{[n]})\) satisfying \(\delta \phi = g - \hat{g}(\emptyset)\) suffices for Theorem~\ref{harris-kleitman}.
Taking instead \(\phi = \big(\sum_{a \in [n]} \sum_{A \subseteq [a - 1]} x^A dx_a\big) * dg\), we morally recover Kleitman's inductive proof \cite{kleitman1966}, where we proceed via the identity
\[
    \langle f, g \rangle = \left\langle \tfrac{\partial}{\partial x_k} f, \tfrac{\partial}{\partial x_k} g \right\rangle + \left\langle f \big|_{x_k = 0}, g \big|_{x_k = 0} \right\rangle.
\]
\end{remark}

\begin{lemma}
\label{path-decomposition}
Let \(\phi \in C^1(3^{[n]})\) be a nonnegative, integer-valued \(1\)-cochain. Then \(\phi\) decomposes into a sum of monotone paths in the discrete cube, each running from a vertex \(A \subseteq [n]\) with \(\delta \phi(A) < 0\) to a vertex \(B \subseteq [n]\) with \(\delta \phi(B) > 0\) such that \(A \subset B\).
\end{lemma}

\begin{proof}
We proceed by induction.
If \(\phi = 0\), the decomposition is empty.
Otherwise, some vertex has nonzero codifferential.
Pick \(A \subseteq [n]\) with \(\delta \phi(A) < 0\); the outflow at \(A\) exceeds the inflow, so by integrality some outgoing edge has weight at least \(1\).
Follow such an edge to a successor; at any intermediate vertex \(B \subseteq [n]\), if \(\delta\phi(B) \le 0\), the inflow is at least balanced by outgoing edges, ensuring we can arrive at some \(B\) such that \(\delta \phi(B) > 0\).
Subtracting this path produces a nonnegative integer \(1\)-cochain of strictly smaller weight, which decomposes by the inductive hypothesis; adjoining the removed path gives the decomposition of \(\phi\).
\end{proof}

One step of this decomposition is illustrated in Fig.~\ref{path-figure}.

\begin{figure}[!b]
\centering
\setlength{\tabcolsep}{3pt}
\begin{tabular}{ccccc}
    \cubeflow{e/v3/1, v1/v12/1, v3/v13/1, v3/v23/1, v13/v123/1}
        {e/{\(-1\)}/below, v1/{\(-1\)}/left, v3/{\(-1\)}/right,
         v12/{\(+1\)}/left, v23/{\(+1\)}/right,
         v123/{\(+1\)}/above}{\(\phi\)} &
    \(=\) &
    \cubeflow{e/v3/1, v3/v13/1, v13/v123/1}
        {e/{\(-1\)}/below, v123/{\(+1\)}/above}{\(\gamma\)} &
    \(+\) &
    \cubeflow{v1/v12/1, v3/v23/1}
        {v1/{\(-1\)}/left, v3/{\(-1\)}/right,
         v12/{\(+1\)}/left, v23/{\(+1\)}/right}{\(\phi'\)}
\end{tabular}
\caption{One step of the decomposition in \(3^{[3]}\), removing a monotone path \(\gamma\) from \(\phi\). Edge labels give the weights of each \(1\)-cochain and vertex labels its codifferential, both omitted when zero.}
\label{path-figure}
\end{figure}

\begin{theorem}[Generalization of {\cite{frankl2023}}]
We first define a \((p, q)\)-matching in a relation \(R \subseteq A \times B\) to be a multiset \(M\) of pairs in \(R\) such that each \(a \in A\) appears in exactly \(p\) pairs of \(M\) and each \(b \in B\) appears in exactly \(q\) pairs of \(M\).

Then for \(\mathcal{F}, \mathcal{G} \subseteq 2^{[n]}\), superset-closed families, there exists a \((p, q)\)-matching between the pairs \((F, G) \in \mathcal{F} \times \mathcal{G}\) with \(F \cup G = [n]\), for \(p = \lvert \mathcal{G} \rvert / \gcd(\lvert \mathcal{F} \rvert, \lvert \mathcal{G} \rvert)\) and \(q = \lvert \mathcal{F} \rvert / \gcd(\lvert \mathcal{F} \rvert, \lvert \mathcal{G} \rvert)\).
\end{theorem}

\begin{proof}
Let \(f = \sum_{F \in \mathcal{F}} \tilde{F}\) and \(g = \sum_{G \in \mathcal{G}} \widetilde{[n] \setminus G}\).

Set \(p = \lvert \mathcal{G} \rvert / \gcd(\lvert \mathcal{F} \rvert, \lvert \mathcal{G} \rvert)\) and \(q = \lvert \mathcal{F} \rvert / \gcd(\lvert \mathcal{F} \rvert, \lvert \mathcal{G} \rvert)\), and let
\[
    h = p f - q g.
\]
Then \(h\) is monotone and integer-valued with \(\hat{h}(\emptyset) = 0\), and by Corollary~\ref{pseudoinverse-nonnegative}, the pseudoinverse \(\delta^\dagger h\) is nonnegative.

Set \(\phi = \delta^\dagger h\).
We modify \(\phi\) to be integer-valued without changing its codifferential or losing nonnegativity.
While the non-integral edges of \(\phi\) contain a directed cycle \(\psi\), set \(\phi\) to be \(\phi + t\psi\), where \(t > 0\) is minimal such that some edge of \(\psi\) becomes integral.
Nonnegativity is preserved, since \(t\) is at most the minimum value of \(\phi\) on edges where \(\psi\) is negative (in the worst case, \(t\) equals that minimum and the corresponding edge becomes zero).
Since \(\delta \psi = 0\), this preserves \(\delta \phi = h\), and the count of non-integral edges strictly decreases, so the process terminates.

If on termination non-integral edges remain, they form a forest.
At a leaf \(A\) of any component, exactly one incident edge has non-integral \(\phi\)-value and the remaining incident edges have integer \(\phi\)-values.
Hence \(\delta \phi(A)\) is the signed sum of an integer and a non-integer, and so non-integer; but \(\delta \phi(A) = h(A) \in \mathbb{Z}\), a contradiction.

Thus \(\phi\) is a nonnegative integer \(1\)-cochain with \(\delta \phi = h\).
By Lemma~\ref{path-decomposition}, \(\phi\) decomposes into a multiset of monotone paths, with each such path corresponding to a pair \((F, G) \in \mathcal{F} \times \mathcal{G}\) with \([n] \setminus G \subseteq F\), equivalently \(F \cup G = [n]\).

Note that these paths alone do not yet form the complete matching.
For \(A \in \mathcal{F}\) with \([n] \setminus A \in \mathcal{G}\) we have \(h(A) = p - q\) rather than \(p\), so such an \(A\) is the endpoint of only \(\max(p - q, 0)\) paths and the starting point of only \(\max(q - p, 0)\), but by adjoining \(\min(p, q)\) copies of the pair \((A, [n] \setminus A)\) for each such \(A\), which satisfies the union condition, we correct both counts to \(p\) and \(q\) respectively.
\end{proof}

\section{Convex Duality}

Given \(f \colon 2^{[n]} \to \mathbb{R}\) and weights \(\alpha \in \mathbb{R}_{\ge 0}^{[n]}\), we consider the optimization problem of finding \(\phi^* \in C^1(3^{[n]})\) that maximizes
\[
    \sum_{k \in [n]} \alpha_k \min(\phi_k)
\]
among all \(\phi \in C^1(3^{[n]})\) subject to the affine constraint \(\delta \phi = f - \hat{f}(\emptyset)\).

To analyze this problem, we recall the Karush-Kuhn-Tucker (KKT) conditions, which give necessary and sufficient conditions for optimality in convex optimization problems with affine equality constraints.

\begin{theorem}[Karush-Kuhn-Tucker conditions {\cite{rockafellar1970}}]
Consider the optimization problem
\[
    \begin{aligned}
        \text{minimize} \quad & f(x) \\
        \text{subject to} \quad & g(x) = 0,
    \end{aligned}
\]
where \(f \colon \mathbb{R}^n \to \mathbb{R}\) is convex (not necessarily differentiable) and \(g \colon \mathbb{R}^n \to \mathbb{R}^m\) is affine.

A point \(x^* \in \mathbb{R}^n\) is a global minimizer if and only if there exists \(\lambda \in \mathbb{R}^m\) such that
\begin{enumerate}
    \item \(g(x^*) = 0\), and
    \item \(0 \in \partial f(x^*) + \nabla g(x^*)^\top \lambda\).
\end{enumerate}
Here \(\partial f\) denotes the subdifferential of \(f\), and \(\nabla g\) denotes the Jacobian of \(g\).
\end{theorem}

Applying the KKT conditions to the optimization problem above, we obtain the following primal, dual, and strong duality statements.

\begin{lemma}
\label{primal}
Given \(f \colon 2^{[n]} \to \mathbb{R}\) and weights \(\alpha \in \mathbb{R}_{\ge 0}^{[n]}\),
then \(\phi^* \in C^1(3^{[n]})\) such that \(\delta \phi^* = f - \hat{f}(\emptyset)\) is a maximizer of
\[
    \sum_{k \in [n]} \alpha_k \min(\phi_k)
\]
for \(\phi \in C^1(3^{[n]})\) subject to the constraint \(\delta \phi = f - \hat{f}(\emptyset)\) if and only if there exists
\[
    \lambda = \sum_{i \in I} c_i \lambda_i
\]
for some \(c_i > 0\) and monotone \(\lambda_i : 2^{[n]} \to \{-1,1\}\) such that:
\begin{enumerate}
    \item for each \(i \in I\), the support of \(d\lambda_i\) is contained in the minima of \(\phi^*\), and
    \item for each \(k \in [n]\), \(\sum_{i \in I} c_i \widehat{\lambda_i}(\{k\}) = \alpha_k\).
\end{enumerate}
\end{lemma}

\begin{proof}
The objective is concave and the constraint affine, and the adjoint of \(\delta\) is \(d\), so by the Karush-Kuhn-Tucker conditions \(\phi^*\) is a maximizer if and only if \(\delta \phi^* = f - \hat{f}(\emptyset)\) and \(d\lambda\) lies in the superdifferential of \(\sum_{k \in [n]} \alpha_k \min(\phi_k)\) at \(\phi^*\), for some \(\lambda \colon 2^{[n]} \to \mathbb{R}\).
Since the superdifferential of \(\min\) is the convex hull of the coordinates attaining it, this says exactly that \(d\lambda\) is nonnegative, is supported on the minima of \(\phi^*\), and satisfies \(\widehat{\lambda}(\{k\}) = \alpha_k\) for each \(k \in [n]\), the last because \((d\lambda)_k\) has mean \(\widehat{\lambda}(\{k\})\).
We can then decompose \(\lambda\) into a sum of monotone functions as required.
\end{proof}

\begin{lemma}
\label{dual}
Given \(f \colon 2^{[n]} \to \mathbb{R}\) and weights \(\alpha \in \mathbb{R}_{\ge 0}^{[n]}\),
then \(\lambda^*\) minimizes \(\langle \lambda, f - \hat{f}(\emptyset) \rangle\) among all \(\lambda \colon 2^{[n]} \to \mathbb{R}\) such that
\[
    \lambda = \sum_{i \in I} c_i \lambda_i
\]
for some \(c_i > 0\) and monotone \(\lambda_i : 2^{[n]} \to \{-1,1\}\) such that for each \(k \in [n]\),
\[
    \sum_{i \in I} c_i \widehat{\lambda_i}(\{k\}) = \alpha_k
\]
if and only if there exists some \(\phi \in C^1(3^{[n]})\) with \(\delta \phi = f - \hat{f}(\emptyset)\) such that for each \(i \in I\), the support of \(d\lambda_i\) is contained in the minima of \(\phi\).
\end{lemma}

\begin{proof}
For any such \(\lambda\) and any \(\phi \in C^1(3^{[n]})\) with \(\delta \phi = f - \hat{f}(\emptyset)\), the adjointness of \(d\) and \(\delta\) gives \begin{align*}
    \langle \lambda, f - \hat{f}(\emptyset) \rangle & = \langle \lambda, \delta \phi \rangle = \langle d\lambda, \phi \rangle = \sum_{k \in [n]} \frac{1}{2^{n - 1}} \sum_{A \subseteq [n] \setminus \{k\}} (d\lambda)_k(A) \phi_k(A) \\
                               & \ge \sum_{k \in [n]} \alpha_k \min(\phi_k),
\end{align*} since each \(d\lambda_i\) is nonnegative and \((d\lambda)_k\) has mean \(\sum_{i \in I} c_i \widehat{\lambda_i}(\{k\}) = \alpha_k\), with equality precisely when \(d\lambda\) is supported on the minima of \(\phi\).
By Lemma~\ref{primal} a maximizer \(\phi^*\) admits such a \(\lambda\), so the two optima agree, and \(\lambda^*\) is a minimizer if and only if it attains this equality against some such \(\phi\).
\end{proof}

\begin{lemma}
\label{strong-duality}
Given \(f \colon 2^{[n]} \to \mathbb{R}\) and weights \(\alpha \in \mathbb{R}_{\ge 0}^{[n]}\),
let \(\phi \in C^1(3^{[n]})\) satisfy \(\delta \phi = f - \hat{f}(\emptyset)\),
and let \(\lambda = \sum_{i \in I} c_i \lambda_i\) for some \(c_i > 0\) and monotone \(\lambda_i \colon 2^{[n]} \to \{-1, 1\}\) with \(\sum_{i \in I} c_i \widehat{\lambda_i}(\{k\}) = \alpha_k\) for each \(k \in [n]\).
Then \[
    \sum_{k \in [n]} \alpha_k \min(\phi_k) = \langle \lambda, f - \hat{f}(\emptyset) \rangle
\] if and only if the support of \(d\lambda\) is contained in the minima of \(\phi\), making \(\phi\) a maximizer and \(\lambda\) a minimizer.
\end{lemma}

\begin{proof}
As in the proof of Lemma~\ref{dual}, we have \[
    \langle \lambda, f - \hat{f}(\emptyset) \rangle = \langle d\lambda, \phi \rangle \ge \sum_{k \in [n]} \alpha_k \min(\phi_k),
\] with equality precisely when \(d\lambda\) is supported on the minima of \(\phi\), and hence each \(d\lambda_i\) is also, being nonnegative.
Thus \(\phi\) is a maximizer by Lemma~\ref{primal} and \(\lambda\) a minimizer by Lemma~\ref{dual}, with their optima equal.
\end{proof}

\begin{remark}
Each maximizer \(\phi \in C^1(3^{[n]})\) for a given weight \(\alpha \in \Delta^{[n]}\) determines a convex region of weights \(A \subseteq \Delta^{[n]}\) for which \(\phi\) is a maximizer.
Moreover, for any such region \(A\), let \(\Phi \subseteq C^1(3^{[n]})\) be the convex region of \(\phi\) that are maximizers at every \(\alpha \in A\), and let \(B\) be the set of all \((\min(\phi_k))_{k \in [n]}\) for all \(\phi \in \Phi\).
Then \(\Delta^{[n]}\) uniquely decomposes into the set of maximal regions \(A\).

Vertices of the decomposition are given by
\[
    \frac{1}{\sum_{k \in [n]} \hat{\lambda}(\{k\})} (\hat{\lambda}(\{k\}))_{k \in [n]},
\]
where \(\lambda \colon 2^{[n]} \to \{-1, 1\}\) is monotone and the support of \(d\lambda\) is contained in the minima of every \(\phi\) that is a maximizer at this vertex.

All non-maximal regions of \(\relint(\Delta^{[n]})\) are given by the intersection of maximal regions containing them, which corresponds to convex combinations of corresponding \(\Phi\) and \(B\).
And finally, for all regions of \(\relint(\Delta^{[n]})\), we have \(\dim(A) + \dim(B) = n - 1\).
\end{remark}

The decomposition of \(\Delta^{[3]}\) for the majority function is given in Fig.~\ref{simplex-figure}.

\begin{figure}[!b]
\centering
\setlength{\tabcolsep}{10pt}
\begin{tabular}{cl}
    \simplexf &
    \(\begin{aligned}
        f                   & = \tfrac{1}{2} x + \tfrac{1}{2} y + \tfrac{1}{2} z - \tfrac{1}{2} xyz \\
        \phi_{\mathrm{i}}   & = \tfrac{1}{2} dx + \tfrac{1}{2} dy + \tfrac{1}{2} dz - \tfrac{1}{2} yz dx \\
        \phi_{\mathrm{ii}}  & = \tfrac{1}{2} dx + \tfrac{1}{2} dy + \tfrac{1}{2} dz - \tfrac{1}{2} xz dy \\
        \phi_{\mathrm{iii}} & = \tfrac{1}{2} dx + \tfrac{1}{2} dy + \tfrac{1}{2} dz - \tfrac{1}{2} xy dz
    \end{aligned}\)
\end{tabular}
\caption{The \(\Delta^{[3]}\) decomposition of \(f\), with representative maximizers given for each maximal region of weights.}
\label{simplex-figure}
\end{figure}

\begin{theorem}
\label{nonnegative}
Let \(f \colon 2^{[n]} \to \mathbb{R}\) be monotone.
Then if \(\phi \in C^1(3^{[n]})\) such that \(\delta \phi = f - \hat{f}(\emptyset)\) is a maximizer for \(\alpha \in \mathbb{R}_{\ge 0}^{[n]}\), we have \(\min(\phi_k) \ge 0\) whenever \(\alpha_k > 0\) for all \(k \in [n]\).
\end{theorem}

\begin{proof}
Fix \(k\) with \(\alpha_k > 0\).
Then by Lemma~\ref{primal}, we have \(\sum_{i \in I} c_i \widehat{\lambda_i}(\{k\}) = \alpha_k\), so some \(\lambda = \lambda_i\) has \(\widehat{\lambda}(\{k\}) > 0\) and thus \((d\lambda)_k \ne 0\), with \(d\lambda\) supported in the minima of \(\phi\).
Then \begin{align*}
      & \sum_{\mathclap{(d\lambda)_k(A) = 1}} f(A \cup \{k\}) - f(A) \\
    = & \sum_{\mathclap{(d(\lambda |_{x_k = 1}))_a(A) = 1}} \phi_a(A \cup \{k\}) \quad - \quad \sum_{\mathclap{(d(\lambda |_{x_k = -1}))_a(A) = 1}} \phi_a(A \cup \{k\}) \quad + \quad \sum_{\mathclap{(d\lambda)_k(A) = 1}} \phi_k(A) \\
      & - \sum_{\mathclap{(d(\lambda |_{x_k = 1}))_a(A) = 1}} \phi_a(A) \quad + \quad \sum_{\mathclap{(d(\lambda |_{x_k = -1}))_a(A) = 1}} \phi_a(A) \quad + \quad \sum_{\mathclap{(d\lambda)_k(A) = 1}} \phi_k(A) \\
    = & \sum_{\mathclap{(d(\lambda |_{x_k = 1}))_a(A) = 1}} \min(\phi_a) - \phi_a(A) \quad + \quad \sum_{\mathclap{(d(\lambda |_{x_k = -1}))_a(A) = 1}} \min(\phi_a) - \phi_a(A \cup \{k\}) \\
      & + \sum_{\mathclap{(d\lambda)_k(A) = 1}} 2 \min(\phi_k).
\end{align*}
Thus we have \( \min(\phi_k) \ge 0 \) by monotonicity.
\end{proof}

\begin{theorem}
\label{containment}
Let \(f, g \colon 2^{[n]} \rightarrow \{-1, 1\}\) be monotone with \(g^{-1}(1) \subseteq f^{-1}(1)\) and \(g\) nonconstant.
Define \(h \colon 2^{[n]} \rightarrow \{-1, 1\}\) such that \(h\) is \(1\) at \([n]\) and \(-1\) otherwise.
Then \[
    \frac{\langle f, g \rangle - \hat{f}(\emptyset) \hat{g}(\emptyset)}{\sum_{k \in [n]} \hat{g}(\{k\})} \ge \frac{\langle f, h \rangle - \hat{f}(\emptyset) \hat{h}(\emptyset)}{\sum_{k \in [n]} \hat{h}(\{k\})}.
\]
\end{theorem}

\begin{proof}
We first calculate,
\[
    \langle f, g \rangle - \hat{f}(\emptyset) \hat{g}(\emptyset) = \langle f - \hat{f}(\emptyset), g + 1 \rangle = \frac{1 - \hat{f}(\emptyset)}{2^{n - 1}} |g^{-1}(1)|.
\]
Each \(A \in g^{-1}(1)\) has at most \(n\) incoming boundary edges, so 
\[
    \sum_{k \in [n]} \hat{g}(\{k\}) \le \frac{n}{2^{n - 1}} |g^{-1}(1)|.
\]
The ratio is thus at least \((1 - \hat{f}(\emptyset))/n\), with equality achieved by \(h\).
\end{proof}

\begin{remark}
Theorem~\ref{containment} follows similarly for \(f, g \colon 2^{[n]} \rightarrow \{-1, 1\}\) monotone with \(g^{-1}(-1) \subseteq f^{-1}(-1)\) and \(g\) nonconstant, with \(h \colon 2^{[n]} \rightarrow \{-1, 1\}\) such that \(h\) is \(-1\) at \(\emptyset\) and \(1\) otherwise.
\end{remark}

\section{Polarity}

\begin{definition}[Polarity]
We say that \(f \colon 2^{[n]} \to \mathbb{R}\) is polar if there exists \(\phi \in C^1(3^{[n]})\) with \(\delta \phi = f - \hat{f}(\emptyset)\) such that every \(\emptyset\)-adjacent edge is minimal, i.e. \(\phi_k(\emptyset) = \min(\phi_k)\) for all \(k \in [n]\).
Likewise, we say \(\phi\) is polar.
\end{definition}

We now give a generalization of Kleitman's conjecture, analogous to Kahn's Conjecture~3.4 in \cite{friedgut2018}, but with an intersection condition in place of antipodality. It similarly implies Kahn's Conjecture~3.7.

\begin{conjecture}
\label{kleitman-cochain}
Let \(f \colon 2^{[n]} \to \{-1, 1\}\) be monotone and intersecting.
Then \(f\) is polar.
\end{conjecture}

An intuitive framing for Conjecture~\ref{kleitman-cochain} comes from the heat flow derivation of \(\delta^\dagger f\): \[
    \delta^\dagger f = \int_{0}^{\infty} d \exp(-t \Delta) f \, dt,
\]
where one might expect, albeit falsely, heat flow to be maximized near \(df\) and minimized near the poles for monotone \(f \colon 2^{[n]} \to \{-1, 1\}\), provided \(df\) is sufficiently smooth.

\begin{theorem}
Conjecture~\ref{kleitman} and Conjecture~\ref{kleitman-cochain} in the case of maximal intersecting \(f \colon 2^{[n]} \to \{-1, 1\}\) are equivalent.
\end{theorem}

\begin{proof}
Assuming Conjecture~\ref{kleitman}, we have \[
    f = \sum_{k \in [n]} \alpha_k x_k + \sum_{A \subset B \subseteq [n]} \beta(A, B) (\tilde{B} - \tilde{A})
\] for some \(\alpha \in \Delta^{[n]}\) and \(\beta \colon \{(A, B) \colon A \subset B \subseteq [n]\} \to \mathbb{R}_{\ge 0}\).
We then have \[
    \phi = \sum_{k \in [n]} \alpha_k dx_k + \sum_{A \subset B \subseteq [n]} \beta(A, B) \gamma(A, B),
\] a \(1\)-cochain satisfying \(\delta \phi = f\), where \(\gamma \colon \{(A, B) \colon A \subset B \subseteq [n]\} \to C^1(3^{[n]})\) gives the \(1\)-cochain of some monotone path from \(A\) to \(B\).
Here \(\phi\) is polar, since its second sum is nonnegative and vanishes at \(\emptyset\), since \(\sum_{k \in [n]} \alpha_k = 1\).

Assuming Conjecture~\ref{kleitman-cochain}, we have \(\delta \phi = f\) for some polar \(\phi \in C^1(3^{[n]})\).
We then have \[
    f = \sum_{k \in [n]} \phi_k(\emptyset) x_k + \sum_{k \in [n]} \sum_{A \subseteq [n] \setminus \{k\}} \left( \phi_k(A) - \phi_k(\emptyset) \right) (\widetilde{A \cup \{k\}} - \tilde{A}),
\] where polarity and Theorem~\ref{nonnegative} ensure nonnegativity, and \(\sum_{k \in [n]} \phi_k(\emptyset) = 1\) as before.
\end{proof}

A quick calculation then gives us the following correlation inequality, from which Chv\'{a}tal's conjecture follows, given Conjecture~\ref{kleitman-cochain} in the maximal intersecting case.
Likewise, so too does Kahn's Conjecture~3.7 follow, though given Conjecture~\ref{kleitman-cochain} in the general case.

\begin{theorem}
\label{correlation}
Let \(f \colon 2^{[n]} \to \mathbb{R}\) be monotone and \(g \colon 2^{[n]} \to \mathbb{R}\) monotone and polar.
Then \[
    \langle f, g \rangle - \hat{f}(\emptyset) \hat{g}(\emptyset) \ge \min_{k \in [n]}\left( \hat{f}(\{k\}) \right) (\hat{g}(\emptyset) - g(\emptyset)).
\]
\end{theorem}

\begin{proof}
\begin{align*}
    \langle f, g \rangle - \hat{f}(\emptyset) \hat{g}(\emptyset) & = \langle df, \phi \rangle \\
                                                                 & \ge \sum_{k \in [n]} \hat{f}(\{k\}) \min(\phi_k) \\
                                                                 & \ge \min_{k \in [n]}\left( \hat{f}(\{k\}) \right) (\hat{g}(\emptyset) - g(\emptyset)).
\end{align*}
\end{proof}

\begin{corollary}
Suppose Conjecture~\ref{kleitman-cochain} holds in the maximal intersecting case.
Then Conjecture~\ref{chvatal} holds.
\end{corollary}

\begin{proof}
Let \(\mathcal{F} \subseteq 2^{[n]}\) be subset-closed, and let \(\mathcal{G} \subseteq \mathcal{F}\) be intersecting.
We then extend \(\mathcal{G}\) to a maximal intersecting \(\mathcal{G}^* \subseteq 2^{[n]}\), and let \(\mathcal{F}^* = \{ [n] \setminus A : A \in \mathcal{F} \}\), which is superset-closed.
Finally, we construct \(f = \sum_{A \in \mathcal{F}^*} \tilde{A}\) and \(g = \sum_{A \in \mathcal{G}^*} \tilde{A}\), both monotone, with \(g\) polar by assumption.

By Theorem~\ref{correlation}, we have \[
    \lvert \mathcal{F} \rvert - \lvert \mathcal{F} \cap \mathcal{G}^* \rvert = \lvert \mathcal{F}^* \cap \mathcal{G}^* \rvert = 2^n \langle f, g \rangle \ge 2^{n - 1} \min_{k \in [n]}\left( \hat{f}(\{k\}) \right) + 2^{n - 1} \hat{f}(\emptyset),
\]
and therefore, \[
    \lvert \mathcal{F} \cap \mathcal{G} \rvert \le \lvert \mathcal{F} \cap \mathcal{G}^* \rvert \le 2^{n - 1} \hat{f}(\emptyset) - 2^{n - 1} \min_{k \in [n]}\left( \hat{f}(\{k\}) \right),
\]
where, writing \(\mathcal{F}_k = \{ A \in \mathcal{F} : k \in A \}\) for the star at \(k\), we have \[
    2^{n - 1} \hat{f}(\emptyset) - 2^{n - 1} \hat{f}(\{k\}) = 2^n \langle f, \frac{1}{2} - \frac{1}{2} x_k \rangle = \lvert \mathcal{F}_k \rvert.
\]
Thus, amongst the largest intersecting subfamilies of \(\mathcal{F}\), we have a star.
\end{proof}

We say that \(\mathcal{F} \subseteq 2^{[n]}\) is a level family if, for some \(k \in [n]\), it consists of all \(A \subseteq [n]\) with \(\lvert A \rvert \ge k\); such a family is intersecting precisely when \(2k > n\).
Here we examine this special case, in which we have a closed form for certain minima.

\begin{theorem}
\label{level}
Let \(f \colon 2^{[n]} \to \{-1, 1\}\) be monotone such that, for some \(k\) with \(2k > n\):
\begin{itemize}
    \item \(f(A) = -1\) for all \(A \subseteq [n]\) with \(\lvert A \rvert < k\), and
    \item \(f(A) = 1\) for all \(A \subseteq [n]\) with \(\lvert A \rvert \ge k\).
\end{itemize}
Then \(f\) is polar, admitting \(\phi = \left( \sum_{a \in [n]} \sum_{A \subseteq [a - 1]} x^A dx_a \right) * df\) with \[
    \min(\phi_a) = \frac{1}{2^{n - a}} \binom{n - a}{k - 1}
\]
for all \(a \in [n]\).
\end{theorem}

\begin{proof}
We construct \[
    \phi = \Big( \sum_{a \in [n]} \sum_{A \subseteq [a - 1]} x^A dx_a \Big) * df.
\]
We thus have \[
    \phi_a = \big(\sum_{A \subseteq [a - 1]} x^A\big) * \frac{\partial}{\partial x_a} f
\] for each \(a \in [n]\), where \(\frac{\partial}{\partial x_a} f\) is the indicator of \(\binom{[n] \setminus \{a\}}{k - 1}\), and \(\sum_{A \subseteq [a - 1]} x^A\) takes the value of \(2^{a - 1}\) on the subcube \(\{ B \subseteq [n] : [a - 1] \subseteq B \}\) and \(0\) elsewhere.
Convolving, the pairs contributing to \(\phi_a(A)\) are exactly those \(C\) with \(\lvert C \rvert = k - 1\) and \(C \cap [a - 1] = A \cap [a - 1]\), the remaining \(k - 1 - \lvert A \cap [a - 1] \rvert\) elements of \(C\) being free among the \(n - a\) coordinates outside \([a]\), giving \[
    \phi_a(A) = \frac{1}{2^{n - a}} \binom{n - a}{k - 1 - \lvert A \cap [a - 1] \rvert}.
\]

Writing \(j = \lvert A \cap [a - 1] \rvert\), which ranges over \(\{0, \dots, a - 1\}\), the sequence \(\binom{n - a}{k - 1 - j}\) is unimodal in \(j\), and so is minimized at \(j = 0\) or \(j = a - 1\).
Now \(\phi_a(\emptyset)\) is the value at \(j = 0\), and \(\binom{N}{r'} \le \binom{N}{r}\) for \(r \le r'\) precisely when \(r + r' \ge N\).
Taking \(r = k - a\), \(r' = k - 1\), and \(N = n - a\), this gives \(\phi_a(\emptyset) = \min(\phi_a)\) exactly when \(2k \ge n + 1\), which holds by assumption.
Hence \(f\) is polar.
\end{proof}

\begin{theorem}
\label{kleitman-restriction}
Let \(f \colon 2^{[n]} \to \{-1, 1\}\) be monotone and nonconstant, and let \(\phi \in C^1(3^A)\) for some \(A \sqcup B = [n]\) such that \(\delta \phi = f \big|_{\{x_b = 0\}_{b \in B}} - \hat{f}(\emptyset)\) and
\[
    \sum_{k \in A} \min(\phi_k) = 1 + \hat{f}(\emptyset).
\]
Then \(f\) is polar.
\end{theorem}

\begin{proof}
We construct
\[
    \psi = \sum_{a \in \{-1, 1\}^A} 1_a \wedge \delta^\dagger \left( f \big|_{\{x_k = a_k\}_{k \in A}} \right),
\] where \(1_a \colon 2^A \to \{0, 1\}\) is the indicator of \(a\).

By direct calculation, we have \begin{align*}
    \delta(\phi + \psi)(C) & = \left( f \big|_{\{x_b = 0\}_{b \in B}} - \hat{f}(\emptyset) \right)(A \cap C) + \delta \delta^\dagger \left( f \big|_{\{x_k = c_k\}_{k \in A}} \right)(B \cap C) \\
                           & = \left( f \big|_{\{x_b = 0\}_{b \in B}}(A \cap C) - \hat{f}(\emptyset) \right) + \left( f(C) - f \big|_{\{x_b = 0\}_{b \in B}}(A \cap C) \right) \\
                           & = f(C) - \hat{f}(\emptyset),
\end{align*}
where \(c \in \{-1, 1\}^{[n]}\) is the vertex corresponding to \(C\).

Moreover, \(\psi\) is nonnegative by Corollary~\ref{pseudoinverse-nonnegative}.
Hence \[
    \sum_{k \in [n]} \min(\phi_k + \psi_k) \ge \sum_{k \in A} \min(\phi_k + \psi_k) \ge 1 + \hat{f}(\emptyset).
\]

For the reverse inequality, let \(\lambda \colon 2^{[n]} \to \{-1, 1\}\) be \(-1\) at \(\emptyset\) and \(1\) elsewhere.
We have \(f(\emptyset) = -1\), since \(f\) is nonconstant, so by Lemma~\ref{dual} for \(\alpha = 1\), we have \[
    \sum_{k \in [n]} \min(\phi_k + \psi_k) \le \langle 2^{n - 1} \lambda, f - \hat{f}(\emptyset) \rangle = \hat{f}(\emptyset) - f(\emptyset) = 1 + \hat{f}(\emptyset).
\]

Thus \(\sum_{k \in [n]} \min(\phi_k + \psi_k) = \langle 2^{n - 1} \lambda, f - \hat{f}(\emptyset) \rangle\), and by Lemma~\ref{strong-duality}, \(f\) is polar.
\end{proof}

\begin{corollary}
\label{two-stars}
Let \(\mathcal{F} \subseteq 2^{[n]}\) be a maximal intersecting family contained in the union of two stars.
Then \(-1 + 2 \sum_{A \in \mathcal{F}} \tilde{A}\) is polar.
\end{corollary}

\begin{proof}
Choose \(a, b \in [n]\) such that \(a \in A\) or \(b \in A\) for all \(A \in \mathcal{F}\).
Let \(f = -1 + 2\sum_{A \in \mathcal{F}} \tilde{A}\) and \(g = f \big|_{\{x_k = 0\}_{k \in [n] \setminus \{a, b\}}}\), both monotone and antipodal since \(\mathcal{F}\) is maximal intersecting.
We thus have \(g(\emptyset) = -1\), \(g(\{a, b\}) = 1\), and \(g(\{a\}) = -g(\{b\})\).

The cochain \(\phi \in C^1(3^{\{a, b\}})\) with \(\phi_a = \frac{1 + g(\{a\})}{2}\) and \(\phi_b = \frac{1 + g(\{b\})}{2}\) satisfies \(\delta \phi = g\), and is minimal on every edge, being constant in each direction.
Finally, \(\sum_{k \in \{a, b\}} \min(\phi_k) = 1 = 1 + \hat{f}(\emptyset)\).
Thus by Theorem~\ref{kleitman-restriction} \(f\) is polar.
\end{proof}

\begin{remark}
Theorem~\ref{kleitman-restriction} provides a natural generalization of the 2-edge and 3-edge cases mentioned by Kleitman in \cite{kleitman1979}.
Regarding Corollary~\ref{two-stars}, Chv\'{a}tal's conjecture in this case was proved by Kleitman and Magnanti~\cite{kleitman1974} and later re-proved by Frankl and Kupavskii~\cite{frankl2023}, but to our knowledge no proof of Kleitman's conjecture in this case has appeared in the literature.
\end{remark}

\begin{theorem}
\label{composition}
Let \(f \colon 2^{[m]} \to \{-1, 1\}\) and \(g_1, \dots, g_m \colon 2^{[n]} \to \{-1, 1\}\) be monotone and polar such that:
\begin{itemize}
    \item \(\widehat{g_i}(\emptyset) = 0\) for all \(i \in [m]\), and
    \item \(\langle g_i, g_j \rangle = 0\) for all \(i \ne j\).
\end{itemize}
Define their composition \(f \circ g \colon 2^{[n]} \to \{-1, 1\}\) by
\[
    (f \circ g)(A) = f(g_1(A), \dots, g_m(A)).
\]
Then \(f \circ g\) is polar.
\end{theorem}

\begin{proof}
By our assumptions on \(g_1, \dots, g_m\) and the equality case of Theorem~\ref{harris-kleitman}, each \(g_i\) depends only on the coordinates of some \(A_i \subseteq [n]\), with \(A_1, \dots, A_m\) pairwise disjoint.
Construct \(\phi \in C^1(3^{[m]})\) with \(\delta \phi = f - \hat{f}(\emptyset)\) satisfying the polarity condition, and for each \(i \in [m]\) take \(\psi_i \in C^1(3^{[n]})\) with \(\delta \psi_i = g_i\) satisfying the polarity condition and \((\psi_i)_k = 0\) for all \(k \notin A_i\).
We then construct
\[
    \omega = \sum_{i \in [m]} (\phi_i \circ g) \wedge \psi_i.
\]
For each \(i \in [m]\) and \(A \subseteq [m] \setminus \{i\}\), the function \(\prod_{j \in A} g_j\) depends only on coordinates outside \(A_i\), on which \(\psi_i\) vanishes, so \(\iota_{\nabla(\prod_{j \in A} g_j)} \psi_i = 0\).
Combined with Theorem~\ref{leibniz}, it follows that
\begin{align*}
    \delta \omega & = \sum_{i \in [m]} \sum_{A \subseteq [m] \setminus \{i\}} \delta \left( \widehat{\phi_i}(A) \left( \prod_{j \in A} g_j \right) \wedge \psi_i \right) \\
                  & = \sum_{i \in [m]} \sum_{A \subseteq [m] \setminus \{i\}} \widehat{\phi_i}(A) \left( \prod_{j \in A} g_j \right) g_i \\
                  & = \sum_{A \subseteq [m]} \sum_{i \in A} \widehat{\phi_i}(A \setminus \{i\}) \prod_{i \in A} g_i \\
                  & = \sum_{\emptyset \subset A \subseteq [m]} \widehat{\phi}(A) \prod_{i \in A} g_i \\
                  & = f \circ g - \hat{f}(\emptyset).
\end{align*}

Moreover, expanding the exterior product gives \[
    (\phi_i \circ g) \wedge \psi_i = \sum_{j \in [n]} \left. (\phi_i \circ g) \, (\psi_i)_j \right|_{x_j = 0} dx_j = \sum_{j \in [n]} (\phi_i \circ g) \, (\psi_i)_j \, dx_j,
\] the restriction being vacuous here.
Each such term is then independently polar, since both \(\phi_i \circ g\) and \((\psi_i)_j\) are nonnegative by Theorem~\ref{nonnegative} and attain their minima at \(\emptyset\).
Hence \(f \circ g\) is polar.
\end{proof}

\begin{definition}[Arbitrability]
\label{arbitrated}
Let \(m \colon 2^{[3]} \to \{-1, 1\}\) denote the median,
\[
    m(x, y, z) = \frac{1}{2} x + \frac{1}{2} y + \frac{1}{2} z - \frac{1}{2} xyz.
\]
We say \(f \colon 2^{[n]} \to \{-1, 1\}\) is arbitrated if:
\begin{enumerate}
    \item \(f = -1\), or \(f = x_k\) for some \(k \in [n]\); or
    \item \(f = m(g_1, g_2, x_k)\) for \(g_1, g_2 \colon 2^{[n]} \to \{-1, 1\}\) arbitrated and \(k \in [n]\) such that \(\langle g_1, x_k \rangle = \langle g_2, x_k \rangle = 0\).
\end{enumerate}
\end{definition}

Note that arbitrated functions are monotone, since \(-1\) and \(x_k\) are and \(m\) preserves monotonicity.

\begin{theorem}
Arbitrated functions are polar.
\end{theorem}

\begin{proof}
We proceed by induction.
The cases \(f = -1\) and \(f = x_k\) are trivial.

For the inductive step, let \(f = m(g_1, g_2, x_k)\) and choose \(\psi_1, \psi_2 \in C^1(3^{[n]})\) satisfying polarity for \(g_1, g_2\) respectively.

We construct \[
    \omega = \frac{1}{2} \psi_1 + \frac{1}{2} \psi_2 + \left( \frac{1}{2} - \frac{1}{2} g_1 g_2 \right) dx_k,
\] for which \begin{align*}
    \delta \omega & = \frac{1}{2} \left( g_1 - \widehat{g_1}(\emptyset) \right) + \frac{1}{2} \left( g_2 - \widehat{g_2}(\emptyset) \right) + \left( \frac{1}{2} - \frac{1}{2} g_1 g_2 \right) x_k \\
                  & = m(g_1, g_2, x_k) - \frac{1}{2} \widehat{g_1}(\emptyset) - \frac{1}{2} \widehat{g_2}(\emptyset) \\
                  & = f - \hat{f}(\emptyset).
\end{align*}
Moreover, \(\frac{1}{2} \psi_1\) and \(\frac{1}{2} \psi_2\) are minimal on their \(\emptyset\)-adjacent edges by the inductive hypothesis, as is \(\left( \frac{1}{2} - \frac{1}{2} g_1 g_2 \right) dx_k\), since every arbitrated function takes the value \(-1\) at \(\emptyset\).
Thus \(f\) is polar.
\end{proof}

\begin{definition}
We say \(f \colon 2^{[n]} \to \{-1, 1\}\) is a nonnegative threshold function if
\[
    f(A) = \operatorname{sign}\left( \sum_{k \in [n]} \alpha_k x_k(A) - \beta \right)
\]
for some \(\alpha \in \mathbb{R}_{\ge 0}^{[n]}\) and \(\beta \in \mathbb{R}_{\ge 0}\), where we have \(\operatorname{sign}(0) = -1\).
\end{definition}

\begin{theorem}
Nonnegative threshold functions are arbitrated.
\end{theorem}

\begin{proof}
Let \(\alpha \in \mathbb{R}_{\ge 0}^{[n]}\) and \(\beta \in \mathbb{R}_{\ge 0}\) define \[
    f(A) = \operatorname{sign}\left( \sum_{k \in [n]} \alpha_k x_k(A) - \beta \right).
\]
We assume \(\alpha_1 \ge \dots \ge \alpha_n\) without loss of generality, and we write \(d = d(\alpha)\) for the largest \(k \in [n]\) such that \(\alpha_k\) is nonzero.
We proceed by induction on \(d\).
The cases \(\alpha = 0\) and \(d = 1\) are trivial, giving \(f = -1\) and \(f = -1\) or \(f = x_1\) respectively.

For the inductive step, suppose \(d \ge 2\) and set \begin{align*}
    \beta_k & = \begin{cases}
        0,                   & \text{if } k = d,     \\
        \alpha_k + \alpha_d, & \text{if } k = d - 1, \\
        \alpha_k,            & \text{otherwise,}
    \end{cases} &
    \gamma_k & = \begin{cases}
        0,                   & \text{if } k = d,     \\
        \alpha_k - \alpha_d, & \text{if } k = d - 1, \\
        \alpha_k,            & \text{otherwise.}
    \end{cases}
\end{align*}
Substituting \(x_d = x_{d - 1}\) and \(x_d = -x_{d - 1}\), we thus have \begin{align*}
    f \big|_{x_d = x_{d - 1}} & = \operatorname{sign}\left( \sum_{k \in [n]} \beta_k x_k - \beta \right), &
    f \big|_{x_d = -x_{d - 1}} & = \operatorname{sign}\left( \sum_{k \in [n]} \gamma_k x_k - \beta \right),
\end{align*}
where \(\gamma_{d - 1} = \alpha_{d - 1} - \alpha_d \ge 0\) since \(\alpha\) is decreasing, so both are again nonnegative threshold functions.
Both \(\beta\) and \(\gamma\) are supported in \([d - 1]\), so both restrictions are arbitrated by the inductive hypothesis, and since neither depends on \(x_d\) we have \(\langle f \big|_{x_d = x_{d - 1}}, x_d \rangle = \langle f \big|_{x_d = -x_{d - 1}}, x_d \rangle = 0\).
Since \[
    f = m\left( f \big|_{x_d = x_{d - 1}}, f \big|_{x_d = -x_{d - 1}}, x_d \right),
\]
we thus have \(f\) arbitrated.
\end{proof}

\section{Concluding Remarks}

The following strengthens Conjecture~\ref{kleitman-cochain}, further requiring all maximal regions \(A \subseteq \Delta^{[n]}\) of the decomposition to touch the center.

\begin{conjecture}
Let \(f \colon 2^{[n]} \to \{-1, 1\}\) be monotone and intersecting, and let \(\alpha \in \mathbb{R}_{>0}^{[n]}\). Then for all maximizers \(\phi \in C^1(3^{[n]})\) of \(\sum_{k \in [n]} \alpha_k \min(\phi_k)\) subject to \(\delta \phi = f - \hat{f}(\emptyset)\), we have every \(\emptyset\)-adjacent edge minimal, i.e. \[
    \phi_k(\emptyset) = \min(\phi_k) \quad \text{for all } k \in [n].
\]
\end{conjecture}

Another possible strengthening of Conjecture~\ref{kleitman-cochain} comes with the introduction of a bias.

\begin{conjecture}
\label{kleitman-biased}
Let \(f \colon 2^{[n]} \to \{-1, 1\}\) be monotone and intersecting.
Then for all \(\alpha \in [0, 1]^{[n]}\), there exist \(\phi_k \colon 2^{[n] \setminus \{k\}} \times 2^{[n] \setminus \{k\}} \to \mathbb{R}\) for each \(k \in [n]\) such that \[
    \left( \delta + \iota_{\nabla \sum_{k \in [n]} \alpha_k x_k} \right) \phi = f - f \big|_{\{x_k = -\alpha_k\}_{k \in [n]}},
\]
where \(\phi = \sum_{k \in [n]} \phi_k \big|_{\{y_j = \alpha_j\}_{j \in [n] \setminus \{k\}}} \, dx_k\) is polar, writing \(y\) for the second argument of \(\phi_k\), so that each restriction leaves a function of the first.
\end{conjecture}

Theorem~\ref{composition} admits a natural variant, further restricting \(f\) to be intersecting while relaxing the balancedness of \(g_1, \dots, g_m\) from \(\widehat{g_i}(\emptyset) = 0\) to \(\widehat{g_i}(\emptyset) \le 0\), which would follow from Conjecture~\ref{kleitman-biased}.

\begin{conjecture}
Let \(f \colon 2^{[m]} \to \{-1, 1\}\) and \(g_1, \dots, g_m \colon 2^{[n]} \to \{-1, 1\}\) be monotone and polar such that:
\begin{itemize}
    \item \(\widehat{g_i}(\emptyset) \le 0\) for all \(i \in [m]\),
    \item \(\langle g_i, g_j \rangle = \widehat{g_i}(\emptyset) \widehat{g_j}(\emptyset)\) for all \(i \ne j\), and
    \item \(f\) is intersecting.
\end{itemize}
Then \(f \circ g\) is polar.
\end{conjecture}

We close with the following conjecture, coming from the heat flow derivation of \(\delta^\dagger\).

\begin{conjecture}
Let \(f \colon 2^{[n]} \to \{-1, 1\}\) be monotone and intersecting. Then there exists \(\phi \in C^1(3^{[n]})\) with \(\delta \phi = f - \hat{f}(\emptyset)\) such that the support of \(df\) is contained in the maxima of \(\phi\); that is, for every \(k \in [n]\) and \(A \subseteq [n] \setminus \{k\}\) with \((df)_k(A) = 1\), we have \(\phi_k(A) = \max(\phi_k)\).
\end{conjecture}

\bibliographystyle{plainnat}
\bibliography{article.bib}

\end{document}